\documentclass[11pt,a4paper]{amsart}
\usepackage[all]{xy}
\usepackage{booktabs,amssymb,longtable}

\makeatletter\@addtoreset{equation}{section} \makeatother

\newtheorem{theorem}[equation]{Theorem}

\newtheorem{lemma}[equation]{Lemma}
\newtheorem{corollary}[equation]{Corollary}
\newtheorem{conjecture}[equation]{Conjecture}

\theoremstyle{definition}
\newtheorem{example}[equation]{Example}
\newtheorem{definition}[equation]{Definition}

\theoremstyle{remark}
\newtheorem{remark}[equation]{Remark}

\author{Dimitra Kosta} \thanks{I would like to thank Ivan Cheltsov for suggesting the
problem to me and for useful comments.}

\title{Factoriality of complete intersections in $\mathbb{P}^5$.}

\begin{document}

 \begin{abstract}
Let $X$ be a complete intersection of two hypersurfaces $F_{n}$ and $F_{k}$ in
$\mathbb{P}^5$ of degree
$n$ and $k$ respectively with $n \geq k$, such that the singularities
of $X$ are nodal and $F_{k}$ is smooth. We prove that if the threefold $X$ has at most
$(n+k-2)(n-1)-1$ singular points, then it is factorial.
\end{abstract}
\maketitle\thispagestyle{empty}%

\section{Introduction}

In this paper we shall extend to the complete intersection setting a recent theorem of
Cheltsov \cite{Factorial_threefold_hypersurfaces}, in which he obtained a sharp bound for
the number
of nodes a threefold hypersurface can have and still be factorial.

Suppose that $X$ is the complete intersection of two hypersurfaces $F_{n}$ and $F_{k}$ in
$\mathbb{P}^5$ of degree
$n$ and $k$ respectively with $n \geq k$, such that $X$ is a nodal threefold. We will
prove the following.

\begin{theorem}
\label{main_theorem}
 Suppose that $F_{k}$ is smooth. Then the threefold $X$ is $\mathbb{Q}$-factorial, when
 $$
 |\text{Sing}(X)| \leq (n+k-2)(n-1)-1 \text{ .}
 $$
\end{theorem}

The next example of a non-factorial nodal complete intersection threefold suggests that
the number of nodes,
that a hypersurface can have while being factorial, should be strictly less than
$(n+k-2)^2$.

\begin{example}
Let $X$ be the complete intersection in $\mathbb{P}^5$ of two smooth hypersurfaces
\begin{eqnarray*}
 F  =  x_3 f_1(x_0,x_1,x_2,x_3,x_4,x_5) + x_4 f_2(x_0,x_1,x_2,x_3,x_4,x_5) + x_5
f_3(x_0,x_1,x_2,x_3,x_4,x_5)  =  0\\
 G  =  x_3 g_1(x_0,x_1,x_2,x_3,x_4,x_5) + x_4 g_2(x_0,x_1,x_2,x_3,x_4,x_5) + x_5
g_3(x_0,x_1,x_2,x_3,x_4,x_5)  =  0
\end{eqnarray*}
 where $f_1,f_2,f_3$ are general hypersurfaces of degree $n-1$ and $g_1,g_2,g_3$ general
hypersurfaces of degree $k-1$.
 Then the singular locus $\text{Sing}(X)$, which is given by the vanishing of the
polynomials
 $$
 x_3 = x_4 = x_5 = f_1g_2 - f_2g_1 = f_1g_3 - f_3g_1 = 0 \text{ ,}
 $$
consists of exactly $(n+k-2)^2$ nodal points and the threefold $X$ is not factorial.
\end{example}

Therefore, we can expect the following stated in \cite{Shokurov_vanishing} to be true.

\begin{conjecture}
 Suppose that $F_{k}$ is smooth. Then the threefold $X$ is $\mathbb{Q}$-factorial, when
 $$
 |\text{Sing}(X)| \leq (n+k-2)(n+k-2)-1 \text{ .}
 $$
\end{conjecture}

The assumption of Theorem~\ref{main_theorem} about the smoothness of $F_k$ is essential,
as Example 28 in \cite{Shokurov_vanishing} suggests.

In the case of a nodal threefold hypersurface in $\mathbb{P}^4$, namely when $k=1$,
several attempts where made towards proving Theorem~\ref{main_theorem}, as one can see in
 \cite{IvanPark} and \cite{Dimitra Kosta}. However, a complete proof for $k=1$ was given
in \cite{Factorial_threefold_hypersurfaces}.

\section{Preliminaries}

        Let $\Sigma$ be a finite subset in $\mathbb{P}^{N}$. The points of
        $\Sigma$ impose independent linear conditions on homogeneous forms
        in $\mathbb{P}^{N}$ of degree $\xi$, if for every point $P$ of the
        set $\Sigma$ there is a homogeneous form on $\mathbb{P}^{N}$ of
        degree $\xi$ that vanishes at every point of the set
        $\Sigma\backslash P$  and does not vanish at the point $P$.

The following result, which relates the notion of $\mathbb{Q}$-factoriality with that of
independent linear conditions,
is due to \cite{Cynk} and was stated in the present form in \cite{Shokurov_vanishing}.
\begin{theorem}
\label{2n+k-6}
 The threefold $X$ is $\mathbb{Q}$-factorial in the case when its singular points impose
independent
 linear conditions on the sections of $H^0(\mathcal{O}_{\mathbb{P}^5}(2n+k-6)|_G)$.
\end{theorem}

        The following result was proved in \cite{Eisenbud-Koh} and follows from a result
of J.Edmonds \cite{edmonds}.

        \begin{theorem}
        \label{theorem:eisenbud-koh}
           The points of $\Sigma$ impose independent linear conditions
        on homogeneous forms of degree $\xi\geq 2$ if at most $\xi k +1$
        points of $\Sigma$ lie in a k-dimensional linear subspace of
        $\mathbb{P}^{N}$.
        \end{theorem}

        By \cite{Bese} and \cite{Davis-Geramita} we also know the following.

        \begin{theorem}
        \label{theorem:bese-davis-geramita}
           Let $\pi:Y\to \mathbb{P}^{2}$ be a blow up of distinct
        points $P_1,...,P_{\delta}$ on $\mathbb{P}^{2}$. Then the linear
        system
        $|\pi^{*}(\mathcal{O}_{\mathbb{P}^{2}}(\xi))-\sum_{i=1}^{\delta}E_{i}|$
        is base-point-free for all $\delta\leq \max{(m(\xi +3-m)-1,m^2)}$,
        where $E_i=\pi^{-1}(P_i)$, $\xi \geq 3$, and $m=\lfloor\frac{\xi +3}{2}\rfloor$,
        if at most $k(\xi +3-k)-2$ points of the set
        ${P_1,P_2,\ldots,P_{\delta}}$ lie on a possibly reducible curve of
        degree $1\leq k\leq m$.
        \end{theorem}

What is next is an application, as stated in \cite{Dimitra Kosta}, of the modern
Cayley-Bacharach theorem
(see \cite{eisenbud-green-harris} or \cite{davis-geramita-orecchia}).

\begin{theorem}
        \label{theorem:cayley-bacharach} Let $\Sigma$ be a subset of
        a zero-dimensional complete intersection of the hypersurfaces $X_1, X_2,..., X_N$
        in $\mathbb{P}^{N}$ of degrees $d_1,...,d_N$ respectively. Then the points of
        $\Sigma$ impose dependent linear conditions on
        homogeneous forms of degree $\sum_{i=1}^{N} deg(X_i)-N-1$ if and only if the
        equality
        $|\Sigma|=\prod_{i=1}^{N}d_i$ holds.
        \end{theorem}

Again due to \cite{Factorial_threefold_hypersurfaces} we have the following.

\begin{theorem}
\label{base_locus}
 Let $\Lambda \subseteq \Sigma$ be a subset, let $\phi: \mathbb{P}^r \dasharrow
\mathbb{P}^m$ be a general projection and let
$$
\mathcal{M} \subset | \mathcal{O}_{\mathbb{P}^n}(t)|
$$
be a linear subsystem that contains all hypersurfaces of degree $t$ that pass through
$\Lambda$. Suppose that
\begin{itemize}
 \item the inequality $|\Lambda| \geq (n+k-2)t + 1$ holds,
 \item the set $\phi(\Lambda)$ is contained in an irreducible reduced curve of degree $t$,
\end{itemize}
where $r>m\geq 2$. Then $\mathcal{M}$ has no base curves and either $m=2$ or $t>n+k-2$.
\end{theorem}

 Finally, next is one of our basic tools, a proof of which can be found in
\cite{Ivan4Dec}.

        \begin{theorem}
        \label{theorem:swapping}
          Let $\Sigma$ be a finite subset in
        $\mathbb{P}^{N}$ that is a disjoint union of finite subsets
        $\Lambda$ and $\Delta$, and $P$ be a point in $\Sigma$. Suppose that
        there is a hypersurface in $\mathbb{P}^{N}$ of degree $\alpha \geq
        1$ that contains all points of the set $\Lambda\backslash P$ and
        does not contain $P$, and for every point $Q$ in the set $\Delta$
        there is a hypersurface in $\mathbb{P}^{N}$ of degree $\beta \geq 1$
        that contains all points of the set $\Sigma\backslash Q$ and does
        not contain the point $Q$. Then there is a hypersurface in
        $\mathbb{P}^{N}$ of degree $\gamma$ that contains the set
        $\Sigma\backslash P$ and does not contain the point $P$, where
        $\gamma$ is a natural number such that $\gamma \geq
        max(\alpha,\beta)$.
        \end{theorem}

\section{Proof of Theorem~\ref{main_theorem}}

Let us consider the complete intersection $X$ of two hypersurfaces $F_n$ and $F_k$ in
$\mathbb{P}^5$
of degrees $n$ and $k$ respectively, with $n \geq k$, such that $X$ is a nodal threefold.
Suppose, furthermore, that $F_k$ is smooth and $X$ has at most $(n+k-2)(n-1)-1$ singular
points. We denote now by $\Sigma \subset \mathbb{P}^5$ the set of singular points of $X$.

\begin{definition}
 We say that the points of a subset $\Gamma \subset \mathbb{P}^r$ have property $\star$
if at most
 $t(n+k-2)$ points of the set $\Gamma$ lie on  a curve in $\mathbb{P}^r$ of degree $t\in
\mathbb{N}$.
\end{definition}

For a proof of the following we refer the reader to \cite{Shokurov_vanishing}.

\begin{lemma}
 The points of the set $\Sigma \subset \mathbb{P}^5$ have property $\star$.
\end{lemma}

According to Theorem~\ref{2n+k-6}, for any point $P \in \Sigma$ we need to prove that
there is a hypersurface of degree $2n+k-6$,
that passes through all the points of the set $\Sigma \backslash P$, but not through the
point $P$.

\begin{remark}
\label{remark}
As we mentioned, the claim of Theorem~\ref{main_theorem} is true, when $k=1$ and thus we
need only consider the case $k \geq 2$.
Furthermore, taking into account the following Lemma, we can assume that $n \geq 5$.
\end{remark}

\begin{lemma}
 The threefold $X$ is $\mathbb{Q}$-factorial, when
 $$
 |\text{Sing}(X)| \leq (n+k-2)(n-1)-1 \text{ and } k \leq n \leq 4 \text{ .}
 $$
\end{lemma}

\begin{proof}
Indeed, we consider the projection $$\psi : \mathbb{P}^5 \dasharrow \Pi \cong
\mathbb{P}^2 \text{ ,}$$ from a general plane $\Gamma$ of $\mathbb{P}^5$ to another
general plane $\Pi \cong \mathbb{P}^2$, that sends the set $\Sigma$
to $\psi (\Sigma)= \Sigma'$. Choose a point $P\in \Sigma$ and put $P' = \psi(P)$. We have
the following cases.
\begin{itemize}
 \item If $2=n \geq k=2$, then $|\Sigma|\leq 1$ and the result holds according to
Theorem~\ref{2n+k-6}.
 \item If $3=n \geq k=2$, then $|\Sigma|\leq 5$ and it imposes independent linear
conditions on forms of degree $2$.
 \item If $3=n \geq k=3$, then $|\Sigma|\leq 7$ and it imposes independent linear
conditions on forms of degree $3$.
 \item If $4=n \geq k=2$, then $|\Sigma|\leq 11$ and at most $4t$ points lie on a curve
in $\mathbb{P}^5$ of degree $t$. So, the 11 points of $\Sigma$ impose independent linear
conditions on forms of degree $4$.

 \item If $4=n \geq k=3$, then $|\Sigma|\leq 14$ and at most $5t$ points lie on a curve
in $\mathbb{P}^5$ of degree $t$.

 If the points of $\Sigma' \subset \Pi$ satisfy property $\star$, then the set
 $\Sigma' \backslash P'$ satisfies the requirements of
Theorem~\ref{theorem:bese-davis-geramita} for $\xi=5$
 and this implies that the set $\Sigma$ imposes independent linear conditions on forms of
degree $5$.

 Suppose on the contrary that the points $\Sigma'$ do not satisfy
Theorem~\ref{theorem:bese-davis-geramita} for $\xi =5$.
In this case there is a curve $C_{2}$ of degree $2$ in $\Pi$ that passes through at
least $11$ points of $\Sigma'$. If we take the cone over $C_{2}$ with vertex $\Gamma$, we
obtain a  hypersurface $f_{2}$ in $\mathbb{P}^5$.
Denote by $\Lambda_2$ the points of $\Sigma$ that lie on $f_2$. From
Theorem~\ref{theorem:cayley-bacharach} it follows that the points of
$\Lambda_2$ impose independent linear conditions on homogeneous forms of degree
$5(2-1)-1=4$, since
$\Lambda_2$ is a subset of the complete intersection of hypersurfaces of degree 2 in
$\mathbb{P}^5$.
The set $|\Sigma \backslash \Lambda_2| \leq 3$ imposes independent linear conditions on
forms of degree 2 and, by applying
Theorem~\ref{theorem:swapping}
to the two disjoint sets $\Lambda_2$ and $\Sigma \backslash \Lambda_2$, we get that the
points of $\Sigma$ impose
independent linear conditions on forms of degree $5$.

 \item $4=n \geq k=4$. Then $|\Sigma|\leq 17$ and at most $6t$ points lie on a curve $C_t
\in \mathbb{P}^5$ of degree $t$.

 If the points of $\Sigma' \subset \Pi$ satisfy property $\star$, then the set
 $\Sigma' \backslash P'$ satisfies the requirements of
Theorem~\ref{theorem:bese-davis-geramita} for $\xi=6$
 and this implies that the set $\Sigma$ imposes independent linear conditions on forms of
degree $6$.

 Suppose on the contrary that the points $\Sigma'$ do not satisfy
Theorem~\ref{theorem:bese-davis-geramita} for $\xi =6$.
In this case there is a curve $C_{2}$ of degree $2$ in $\Pi$ that passes through at
least $13$ points of $\Sigma'$. If we take the cone over $C_{2}$ with vertex $\Gamma$, we
obtain a  hypersurface $f_{2}$ in $\mathbb{P}^5$.
Denote by $\Lambda_2$ the points of $\Sigma$ that lie on $f_2$. From
Theorem~\ref{theorem:cayley-bacharach} it follows that the points of
$\Lambda_2$ impose independent linear conditions on homogeneous forms of degree
$5(2-1)-1=4$, since
$\Lambda_2$ is a subset of the complete intersection of hypersurfaces of degree 2 in
$\mathbb{P}^5$.
The set $|\Sigma \backslash \Lambda_2| \leq 4$ imposes independent linear conditions on
forms of degree 2 and, by applying
Theorem~\ref{theorem:swapping}
to the two disjoint sets $\Lambda_2$ and $\Sigma \backslash \Lambda_2$, we get that the
points of $\Sigma$ impose
independent linear conditions on forms of degree $6$.
\end{itemize}
As we saw above, for $3 \leq n \leq 5$ the points of $\Sigma$ impose independent linear
conditions on forms of degree $2n+k-6$, and thus, by Theorem~\ref{2n+k-6}, the threefold
$X$ is $\mathbb{Q}$-factorial.
\end{proof}

\begin{lemma}
\label{sing_plane}
 Suppose that all the singularities of $X$ lie on a plane $\Pi \subset \mathbb{P}^5$.
 Then for any point $P \in \Sigma$ there is hypersurface of degree $(2n+k-6)$ that
contains $\Sigma \backslash P$,
 but does not contain the point $P$.
\end{lemma}

\begin{proof}
By Remark~\ref{remark}, we can see that $\xi=2n+k-6 \geq 6$. Also, we have
$$
|\Sigma \backslash P| \leq \text{max} \left\lbrace
 \lfloor \frac{2n+k-3}{2} \rfloor ( 2n+k-3-\lfloor \frac{2n+k-3}{2} \rfloor)-1,
 \lfloor \frac{2n+k-3}{2} \rfloor^2
 \right\rbrace \text{ ,}
$$ for $k\geq2$ and $n\geq5$. In order to show that at most
$t(2n+k-3-t)-2$ points of $\Sigma$ lie on a curve of degree $t$ in $\Pi$, it is enough to
show that
$$
t(2n+k-3-t)-2 \geq t(n+k-2) \Longleftrightarrow  t(n-t-1) \geq 2  \text{,  for all  } t
\leq \frac{2n+k-3}{2} \text{ .}
$$
For $t=1$ the inequality holds, since $n \geq 5$, and we can assume that $t \geq 2$. It
remains to show that
$t < n-1$. Suppose on the contrary that $t \geq n-1$. The quantity $t(2n+k-3-t)-2$ rises
for all
$n-1 \leq t \leq \lfloor \frac{2n+k-3}{2} \rfloor$ and we have
$$
|\Sigma \backslash P| \leq (n-1)(n+k-2)-2 \leq t(2n+k-3-t)-2 \text{ .}
$$
Therefore we see that the requirement of Theorem~\ref{theorem:bese-davis-geramita}, that
at most $t(2n+k-3-t)-2$
points of $\Sigma$ lie on a curve of degree $t$ in $\Pi$ is satisfied by the set $\Sigma
\backslash P$ for all $t\leq \frac{2n+k-3}{2}$ .
So there is a hypersurface of degree $(2n+k-6)$ that contains $\Sigma \backslash P$,
but does not contain point $P$.
\end{proof}

Taking into account Theorem~\ref{base_locus}, we can reduce to the case $\Sigma$
is a finite set in $\mathbb{P}^3$, such that at most $(n+k-2)t$ of its points
are contained in a curve in $\mathbb{P}^3$ of degree $t \in \mathbb{N}$. Now fix a
general plane $\Pi \in \mathbb{P}^3$ and let
$$
\phi: \mathbb{P}^3 \dasharrow \Pi \cong \mathbb{P}^2
$$
be a projection from a sufficiently general point $O \in \mathbb{P}^3$. Denote by
$\Sigma' = \phi(\Sigma)$ and $P'= \phi(P)$.

\begin{lemma}
 Suppose that the points of $\Sigma' \subseteq \Pi$ have the property $\star$. Then there
is a hypersurface of degree $2n+k-6$
 that contains $\Sigma \backslash P$ and does not contain $P$.
\end{lemma}

\begin{proof}
The points of the set $\Sigma'$ satisfy the requirements of
Theorem~\ref{theorem:bese-davis-geramita},
following the proof of Lemma~\ref{sing_plane}. Thus, there is a curve $C$ in $\Pi$ of
degree $2n+k-6$, that passes through
all the points of the the set $\Sigma' \backslash P'$, but not through the point $P'$. By
taking the cone in $\mathbb{P}^3$
over the curve $C$ with vertex $O$, we obtain the required hypersurface.
\end{proof}

We may assume then, that the points of the set $\Sigma' \subseteq \Pi$ do not have
property $\star$. Then there is a subset
$\Lambda_r^1 \subseteq \Sigma$ with $|\Lambda_r^1| > r(n+k-2)$, but after projection the
points
$$
\phi(\Lambda_r^1) \subseteq \Sigma' \subset \Pi \cong \mathbb{P}^2
$$
are contained in a curve $C_r \subseteq \Pi$ of degree $r$. Moreover, we may assume that
$r$ is the smallest
natural number, such that at least $(n+k-2)r +1$ points of $\Sigma'$ lie on a curve of
degree $r$, which implies that the curve $C_r$ is irreducible and reduced.

By repeating how we constructed $\Lambda_r^1$, we obtain a non-empty disjoint union of
subsets
$$
\Lambda = \bigcup_{j=r}^{l} \bigcup_{i=1}^{c_j} \Lambda_j^i \subseteq \Sigma \text{ ,}
$$
such that $|\Lambda_j^i|>j(n+k-2)$, the points of the set
$$
\phi(\Lambda_j^i) \subseteq \Sigma'
$$
are contained in an irreducible curve in $\Pi$ of degree $j$, and the points of the subset
$$
\phi(\Sigma \backslash \Lambda) \subsetneq \Sigma' \subset \Pi \cong \mathbb{P}^2
$$
have property $\star$, where $c_j \geq 0$.
Let $\Xi_j^i$ be the base locus of the linear subsystem in
$|\mathcal{O}_{\mathbb{P}^3}(j)|$ of all
surfaces of degree $j$ passing through the set $\Lambda_j^i$. Then according to
Theorem~\ref{base_locus}, the base locus $\Xi_j^i$
is a finite set of points and we have $c_r >0$ and
$$
|\Sigma \backslash \Lambda| < (n-1)(n+k-2)- \sum_{i=r}^{l} i(n+k-2)c_i= (n+k-2) \left(
n-1 -\sum_{i=r}^{l} ic_i \right) \text{ .}
$$
\begin{corollary}
\label{bound-sum}
 The inequality $\sum_{i=r}^{l} ic_i \leq n-2$ holds.
\end{corollary}

Put $\Delta = \Sigma \cap (\cup_{j=r}^{l} \cup_{i=1}^{c_j} \Xi_j^i)$. Then $\Lambda
\subseteq  \Delta \subseteq \Sigma$.

\begin{lemma}
\label{delta}
 The points of the set $\Delta$ impose independent linear conditions on forms of degree
$2n+k-6$.
\end{lemma}

\begin{proof}
 We have the exact sequence
 $$
 0 \longrightarrow \mathcal{I}_{\Delta} \otimes \mathcal{O}_{\mathbb{P}^3}(2n+k-6)
 \longrightarrow \mathcal{O}_{\mathbb{P}^3}(2n+k-6) \longrightarrow \mathcal{O}_{\Delta}
 \longrightarrow 0 \text{ ,}
 $$
 where $\mathcal{I}_{\Delta}$ is the ideal sheaf of the closed subscheme $\Delta$ of
$\mathbb{P}^3$.
 Then the points of $\Delta$ impose independent linear conditions on forms of degree
$2n+k-6$,
 if and only if
 $$
 h^1 \left( \mathcal{I}_{\Delta} \otimes \mathcal{O}_{\mathbb{P}^3}(2n+k-6) \right) =0
\text{ .}
 $$
 We assume on the contrary that $h^1 \left( \mathcal{I}_{\Delta} \otimes
\mathcal{O}_{\mathbb{P}^3}(2n+k-6) \right) \neq 0$.
 Let $\mathcal{M}$ be a linear subsystem in $|\mathcal{O}_{\mathbb{P}^3}(n-2)|$ that
contains all surfaces that pass through
 all points of the set $\Delta$. Then the base locus of $\mathcal{M}$ is
zero-dimensional, since
 $\sum_{i=r}^{l} ic_i \leq n-2$ and
 $$
 \Delta \subseteq \cup_{j=r}^{l} \cup_{i=1}^{c_j} \Xi_j^i \text{ ,}
 $$
 but $\Xi_j^i$ is a zero-dimensional base locus of a linear subsystem of
$|\mathcal{O}_{\mathbb{P}^3}(j) |$.
 Let $\Gamma$ be the complete intersection
 $$
 \Gamma = M_1 \cdot M_2 \cdot M_3 \text{ ,}
 $$
 of three general surfaces $M_1, M_2, M_3$ in $\mathcal{M}$.
 Then $\Gamma$ is zero-dimensional and $\Delta$ is closed subscheme of $\Gamma$.
 Let
 $$\mathcal{I}_{\Upsilon}= \text{Ann} \left( \mathcal{I}_{\Delta} / \mathcal{I}_{\Gamma}
\right) \text{ .} $$
 Then
 $$
 0 \neq h^1 \left( \mathcal{I}_{\Delta} \otimes \mathcal{O}_{\mathbb{P}^3}(2n+k-6) \right)
 = h^0 \left( \mathcal{I}_{\Upsilon} \otimes \mathcal{O}_{\mathbb{P}^3}(n-k-4) \right)
 - h^0 \left( \mathcal{I}_{\Gamma} \otimes \mathcal{O}_{\mathbb{P}^3}(n-k-4) \right)
\text{ .}
 $$
 Therefore $h^0 \left( \mathcal{I}_{\Gamma} \otimes \mathcal{O}_{\mathbb{P}^3}(n-k-4)
\right) \neq 0$ and there is a surface
 $F \in |\mathcal{I}_{\Upsilon} \otimes \mathcal{O}_{\mathbb{P}^3}(n-k-4)|$. We have
 $$
 (n-k-4)(n-2)^2 = F \cdot M_2 \cdot M_3 \geq h^0(\mathcal{O}_{\Upsilon}) =
 h^0(\mathcal{O}_{\Gamma})-h^0(\mathcal{O}_{\Delta}) = (n-2)^3-|\Delta| \text{ ,}
 $$
 which implies $|\Delta| \geq (k+2)(n-2)^2$. But $|\Delta| \leq |\Sigma| < (n-1)(n+k-2)$,
which is impossible since $k \geq 2$ and $n \geq 5$.
\end{proof}

We see that $\Delta \subsetneq \Sigma$. Put $\Gamma= \Sigma \backslash \Delta$ and
$d=2n+k-6 - \sum_{i=r}^{l} ic_i$.

\begin{lemma}
\label{d>3}
 The inequality $d \geq 3$ holds.
\end{lemma}

\begin{proof}
 Suppose that $d \leq 2$. Since $\sum_{i=r}^{l} i c_i \leq n-2$ due to
Corollary~\ref{bound-sum}, we have
 $$
 2 \geq d = 2n+k-6-\sum_{i=r}^{l} i c_i \geq 2n+k-6 - (n-2) = n+k-4 \geq 3 \text{ ,}
 $$
 which is impossible.
\end{proof}

For the number of points of the set $\Gamma'$ we have
$$
|\Gamma'| = |\Gamma| \leq |\Sigma \backslash \Lambda| \leq (n+k-2) \left(
n-1-\sum_{i=r}^{l} ic_i \right) - 2 \text{ ,}
$$
and for $d=2n+k-6 - \sum_{i=r}^{l}ic_i$, since $n \geq 5$ and $k \geq 2$, we get
$$
|\Gamma'| \leq (n+k-2) \left( n-1-\sum_{i=r}^{l} ic_i \right) - 2 \leq
\text{max} \left\lbrace
\lfloor \frac{d+3}{2} \rfloor \left( d+3- \lfloor \frac{d+3}{2} \rfloor \right)-1,
\lfloor \frac{d+3}{2} \rfloor^2
\right\rbrace \text{ .}
$$

\begin{lemma}
\label{d-line}
 If the points of the set $\Gamma$ impose dependent linear conditions on forms of degree
$d$, then
 at most $d$ points of the set $\Gamma'$ lie on a line in $\Pi \cong \mathbb{P}^2$.
\end{lemma}

\begin{proof}
 Let us assume on the contrary that there is a line that contains at least $d+1$ points
of $\Gamma$.
 Since the points of $\Gamma$ satisfy property $\star$, at most $n+k-2$ of its points lie
on a line, thus
 $$
 n+k-2 \geq d+1 = 2n+k-6 - \sum_{i=r}^{l}ic_i +1 \text{ ,}
 $$
 which along with Corollary~\ref{bound-sum} implies that
 $$
 n-3 \leq \sum_{i=r}^{l}ic_i \leq n-2 \text{ .}
 $$
 If $\sum_{i=r}^{l}ic_i = n-2$, then $|\Gamma| \leq n+k-4$ and we get a contradiction
 as no more than $n+k-4 < d+1$ points can lie on a line.
 If $\sum_{i=r}^{l}ic_i = n-3$, then $|\Gamma| \leq 2(n+k-3)$ and according to
Theorem~\ref{theorem:eisenbud-koh}
 the points of $\Gamma$ impose independent linear conditions on forms of degree
$d=n+k-3$, which contradicts our assumption.
 By Theorem~\ref{base_locus} the number of points of $\Gamma'$ that can lie on a line
$\Pi \cong \mathbb{P}^2$ is at most $d$.
\end{proof}

\begin{lemma}
\label{<t(d+3-t)-2}
 At most $$t(d+3-t)-2$$ points of the set $\Gamma'$ lie on a curve in $\Pi \cong
\mathbb{P}^2$
 of degree $t$, for every $t \leq \frac{d+3}{2}$.
\end{lemma}

\begin{proof}
 We need to check the condition that at most $t(d+3-t)-2$ points of $\Gamma'$ lie on a
curve of degree $t$
 only for $2 \leq t \leq \frac{d+3}{2}$, such that
 $$
 t(d+3-t)-2 < |\Gamma'| \text{ .}
 $$
 Because the set $\Gamma'$ satisfies property $\star$, at most $(n+k-2)t$ of its points
can lie on a curve of degree $t$
 and therefore it is enough to prove that
 $$
 t(d+3-t)-2 \geq (n+k-2)t  \Longleftrightarrow t\left(  n-1-\sum_{i=r}^{l}ic_i - t
\right)  \geq 2
 \text{ ,  for all  } 2 \leq t \leq \frac{d+3}{2} \text{ .}
 $$
 As we saw Lemma~\ref{d-line} implies that $t \geq 2$ and we only need to show that $t <
n-1- \sum_{i=r}^{l}ic_i$.
 Suppose that
 $$n-1- \sum_{i=r}^{l}ic_i \leq t \leq \frac{d+3}{2} \text{ ,}$$
 then
 $$
 ( n-1- \sum_{i=r}^{l}ic_i ) (n+k-2) = ( n-1- \sum_{i=r}^{l}ic_i ) ( d + 3 -
(n-1-\sum_{i=r}^{l}ic_i ) ) - 2 \leq t(d+3-t)-2 \text{ ,}
 $$
 since the quantity $t(d+3-t)-2$ increases, as $t \leq \frac{d+3}{2}$ increases.
 But then
 $$
 ( n-1- \sum_{i=r}^{l}ic_i ) (n+k-2) -2\leq t(d+3-t)-2 < |\Gamma'|
 \leq ( n-1- \sum_{i=r}^{l}ic_i )(n+k-2)-2 \text{ ,}
 $$
 which is a contradiction.
\end{proof}

\begin{lemma}
 The points of the set $\Sigma$ impose independent linear conditions on homogeneous forms
of degree $2n+k-6$.
\end{lemma}

\begin{proof}
 According to Lemma~\ref{d>3} and Lemma~\ref{<t(d+3-t)-2} all the requirements of
Theorem~\ref{theorem:bese-davis-geramita}
 for $\xi = d$ are satisfied and thus, the points of $\Gamma$ impose independent linear
conditions on homogeneous forms of degree $d$.
 Hence, for any point $Q$ in $\Gamma$, there is a hypersurface $G_{Q}$ of degree $d$,
 such that $G_{Q}(\Gamma \backslash Q)=0$ and $G_{Q}(Q) \neq 0$.

 Furthermore, by the way the set $\Delta$ was constructed, there is a form $F$ of degree
$\sum_{i=r}^{l}ic_i$ in $\mathbb{P}^3$,
 that vanishes at every point of the set  $\Delta$, but does not vanish at any point of
the set $\Gamma$.

 Therefore, for any point $Q \in \Gamma$ we obtain a hypersurface $FG_{Q}$  of degree
$2n+k-6$, such that
 $$
 FG_{Q}(\Sigma) =0 \text{  and  } FG_{Q}(Q) \neq 0 \text{ .}
 $$

Also, by Lemma~\ref{delta}, for any point $R \in \Delta$  there is a hypersurface of
degree $2n+k-6$ that passes through all points of $\Delta\backslash R$, except for the
point $R$.

By applying
Theorem~\ref{theorem:swapping}
to the two disjoint sets $\Delta$ and $\Gamma$, we prove the Lemma.
\end{proof}

        \vskip 1cm

        \noindent
           School of Mathematics\\%
           The University of Edinburgh\\%
           Kings Buildings, Mayfield Road\\%
           Edinburgh EH9 3JZ, UK\\%
           \vskip 0.2cm

           \noindent\texttt{D.Kosta@sms.ed.ac.uk}

\end{document}